\documentclass[12pt]{proc-l}

\usepackage{latexsym}
\usepackage{amsthm}
\usepackage{amsmath}
\usepackage{amssymb}
\usepackage{ulem}
\usepackage{color}

\newtheorem*{theorem*}{Theorem}
\newtheorem*{lemma*}{Lemma}
\newtheorem*{claim*}{Claim}
\newtheorem*{exercise*}{Exercise}
\newtheorem*{note*}{Note}
\newtheorem*{example*}{Example}
\newtheorem*{problem*}{Problem}
\newtheorem*{solution*}{Solution}
\newtheorem*{remark*}{Remark}

\newtheorem{theorem}{Theorem} 
\newtheorem{corollary}[theorem]{Corollary}
\newtheorem{remark}[theorem]{Remark}

\newtheorem{lemma}[theorem]{Lemma}

\newtheorem{proposition}[theorem]{Proposition}

\numberwithin{theorem}{section} \numberwithin{equation}{section}

\newcommand{\ssum}[2] {\overset{#2}{\underset{#1}{\sum}}}  
\newcommand{\nature}{\ensuremath{\mathbb{N}}}   
\newcommand{\integer}{\ensuremath{\mathbb{Z}}}  
\newcommand{\Z}{\integer}

\newcommand{\Q}{\mathbb{Q}}

\newcommand{\co}{s}

\newcommand{\wout}{\backslash}                          

\begin{document}
\title[The triangular theorem of eight]{The triangular theorem of eight and \\
representation by quadratic polynomials}
\author{Wieb Bosma}
\email{bosma@math.ru.nl}
\address{Radboud Universiteit, Heijendaalseweg 135, 6525 AJ Nijmegen, The Netherlands}
\author{Ben Kane}
\email{bkane@mi.uni-koeln.de}
\address{Radboud Universiteit, Heijendaalseweg 135, 6525 AJ Nijmegen, The Netherlands}
\curraddr{Mathematical Institute, University of Cologne, Weyertal 86-90, 50931 Cologne, Germany}
\date{\today}
\begin{abstract}
We investigate here the representability of integers as sums of triangular numbers, where the $n$-th triangular number is given by $T_n=n(n+1)/2$. In particular, we show that $f(x_1, x_2, \ldots, x_k)=b_1T_{x_1}+\cdots+b_kT_{x_k}$, 
for fixed positive integers $b_1, b_2, \ldots, b_k$, represents every
nonnegative integer if and only if it represents $1$, $2$, $4$, $5$, and $8$.
Moreover, if `cross-terms' are allowed in $f$, we show that no finite set
of positive integers can play an analogous role, in turn showing that there is no overarching finiteness theorem which generalizes the statement from positive definite quadratic forms to totally positive quadratic polynomials.
\end{abstract}
\maketitle

\section{Introduction}\label{introsection}
\noindent
In 1638 Fermat claimed that every number is a sum of at most three triangular numbers, four square numbers, and in general $k$ polygonal numbers of order $k$.  The $n$-th
polygonal number of order $k$ is $\frac{(k-2)n^2-(k-4)n}{2}$, so the $n$-th
triangular number is $T_n:=\frac{n(n+1)}{2}$, where we include $T_0=0$ for simplicity.  
The claim for four squares was shown by Lagrange. 
\begin{theorem*}[Lagrange, 1770]
Every positive integer is the sum of four squares.
\end{theorem*}
\noindent
Gauss wrote ``$\text{Eureka, }\triangle+\triangle+\triangle=n$'' in his mathematical diary on July 10, 1796.  
\begin{theorem*}[Gauss, 1796]
Every positive integer is the sum of three triangular numbers.
\end{theorem*}
\noindent
The first proof of the full assertion of Fermat was given by Cauchy in 1813
\cite{Cauchy1}, cf.\ \cite{Nathanson1}.  

For a more complete history of related questions about sums of figurate 
numbers and some new results, see Duke's survey paper \cite{Duke3}.  

\medskip\noindent
The current paper concerns questions of representability of integers by 
quadratic polynomials.
If $f=f(x)=f(x_1, x_2, \ldots,x_k)$ is a rational polynomial in $k$ variables,
it {\it represents} the integer $n$ if there exist
integers $n_i$ such that $n=f(n_1, n_2, \ldots, n_k)$,
and it {\it oddly represents} the integer $n$
if there exist odd integers $n_i$ such that $f(n_1,n_2,\ldots,n_k)=n$.
If $f$ represents every element of a set $\mathcal{Z}$ of integers,
it is said to {\it represent} $\mathcal{Z}$.

If we let $S=S_x$ be the square polynomial $x^2$,
and let $T=T_{x}$ denote the triangular polynomial $(x^2+x)/2$,
the theorems of Lagrange and Gauss state that the positive
integers are represented by $S_w+S_x+S_y+S_z$, and by $T_x+T_y+T_z$.

In 1917, Ramanujan extended the question about four squares to ask
for which choices of quadruples $b=(b_1,b_2,b_3,b_4)$ of integers
the form $b_1S_w+b_2S_x+b_3S_y+b_4S_z$ represents every positive integer;
we shall refer to these as \begin{it}universal diagonal forms\end{it}.
He gave a list of 55 choices of $b$ which he claimed to be the complete
list of universal quarternary diagonal forms; 54 of them turned out to be 
universal and this list is complete, as proven by Dickson \cite{Dickson}.  

Recently, Conway and Schneeberger proved in unpublished work a nice
classification for universal positive definite quadratic forms whose corresponding matrices have integer entries. This answers the question of representability
by positive definite homogeneous quadratic polynomials with {\it even}
off-diagonal coefficients.
\begin{theorem*}[Conway-Schneeberger]
A positive definite quadratic form $Q(x)=x^tAx$, where $A$ is a positive symmetric matrix with integer coefficients, represents every positive integer if and only if it represents the integers $1,2,3,5,6, 7, 10, 14,$ and $15$.  
\end{theorem*}
\noindent
Bhargava gave a simpler proof of the Conway-Schneeberger 15-Theorem in \cite{Bhargava1}, and showed more generally that representability of any $\mathcal{Z}$ by
such form can always be checked
on a finite subset $\mathcal{Y}$. In addition, he exhibited $\mathcal{Y}$ for $\mathcal{Z}$ consisting of all
odd integers and for $\mathcal{Z}$ consisting of all primes.  

More recently, Bhargava and Hanke \cite{BhargavaHanke1} have shown the 290-Theorem, providing the necessary set (the largest
element of which is 290) for universal forms when the corresponding matrix
is half integral, that is, for totally positive integer quadratic forms.

\medskip\noindent
In 1863, Liouville \cite{Liouville1} proved the following generalization of 
Gauss's theorem,
similar to Ramanujan's generalization of Lagrange's Four Squares Theorem.
\begin{theorem*}[Liouville]
Let $a,b,c$ be positive integers with $a\leq b\leq c$.  Then every positive
integer is represented by $aT_{x}+b T_{y}+cT_{z}$ if and only if $(a,b,c)$
is one of the following: 
$$
(1,1,1),\, (1,1,2),\, (1,1,4),\, (1,1,5),\, (1,2,2),\, (1,2,3),\, (1,2,4).
$$
\end{theorem*}
\noindent
We will first prove a finiteness theorem similar to the results of the Conway-Schneeberger 15-Theorem or the Bhargava-Hanke 290-Theorem for sums of triangular numbers. 

\begin{theorem}\label{thm8}
If $b_1,\dots, b_k$ is a sequence of positive integers then
$\ssum{i=1}{k} b_iT_{x_i}$ represents every nonnegative integer if and only if it represents $1$, $2$, $4$, $5$, and $8$.
\end{theorem}
\noindent
Since $8T_x=(2x+1)^2-1$, clearly
$\ssum{i=1}{k} b_iT_{n_i}=n$ if and only if
$\ssum{i=1}{k}b_i(2n_i+1)^2=8n+\ssum{i=1}{k} b_i$.
Hence there is a close correspondence between representability
by triangular polynomials and odd representability by diagonal quadratic
forms.
\begin{corollary}\label{cor8}
If $b_1,\dots, b_k$ is a sequence of positive integers with sum $B$, then
$\ssum{i=1}{k}b_ix_i^2$ oddly represents every integer of the form $8n+B$
with $n\geq 0$ if and only if it oddly represents $8+B$, $16+B$, $32+B$, $40+B$, and $64+B$. 
\end{corollary}
\noindent
It is not so difficult to establish Theorem \ref{thm8} with the escalator
techniques of Bhargava (and Liouville). We will prove a stronger statement
in Section 2:
if the integers $1,2,4,5,$ and $8$ are represented by the triangular form, then $n$ is represented very many times unless $n+1$ has high 3-divisibility.  










\medskip\noindent
We now turn to more general quadratic polynomials.
Let $f$ be a quadratic polynomial in $\Q[x_1,x_2,\dots, x_k]$; 
then $f$ is a 
\begin{it}normalized totally positive\end{it} quadratic polynomial if 
the image of $\integer^k$ under $f$ consists of non-negative integers, while $f(x)=0$ for some $x\in \integer^k$. 
Note that clearly $S_x=x^2$ is normalized totally positive, as is $T_x$:
$T_0=0, T_1=1, T_2=3$ are the first of the increasing sequence
of triangular numbers, and $T_{-m}=T_{m-1}$ for positive $m$.

It turns out that no finiteness theorem will hold in general
for normalized totally positive quadratic polynomials, and moreover that 
checking no proper subset will suffice.
\begin{proposition}\label{quadpolythm}
Let $\mathcal{Z}$ be a subset of the positive integers. For every proper subset 
$\mathcal{Y}\subsetneq \mathcal{Z}$ there exists a normalized totally positive quadratic 
polynomial that represents $\mathcal{Y}$ but does not represent $\mathcal{Z}$.
\end{proposition}
\noindent
Proposition \ref{quadpolythm} will follow directly from the corresponding 
result for {\it triangular sums with cross terms}.
This class corresponds to integral quadratic forms with even off-diagonal terms, just
as the ordinary triangular sums correspond to diagonal quadratic forms.
We refer to Section 3 for a precise definition of this subclass of
quadratic polynomials.

In Section 4 we construct a `norm' $m$  on this class that
restores finite representability.
\begin{theorem}\label{CrossmBoundThm}
Fix an integer $m$ and a subset $\mathcal{Z}$ of the positive integers.  Then there is 
a finite subset $\mathcal{Y}_{m}\subset \mathcal{Z}$, depending only on $m$ and $\mathcal{Z}$, such that 
every triangular sum $t$ with cross terms satisfying $m(t)\leq m$ 
represents $\mathcal{Z}$ if and only if it represents $\mathcal{Y}_{m}$.

Moreover, for $\mathcal{Z}$ equal to the positive integers, we find that $\max \mathcal{Y}_{m}\gg m^2$.
\end{theorem}
\noindent
It may be of interest to investigate the growth of $\max \mathcal{Y}_m$; see Remark \ref{Crossmrmk}. 

\section{Theorem of Eight}\label{Thm8Section}
\noindent
For background information on quadratic forms and genus theory, a good source is \cite{Jones1}.  We prove 
Theorem \ref{repsizethm}, by using a standard argument 
to show that the theorem is equivalent to a statement about (diagonal) 
quadratic forms, and then prove the corresponding result for quadratic forms.  
We will only need some elementary results about quadratic forms and a 
theorem of Siegel to show the desired result.
Theorem \ref{thm8} and Corollary \ref{cor8} follow immediately.

We will first introduce some useful notation and definitions.  
We abbreviate $t(x)=t(x_1, x_2, \ldots, x_k):=\sum b_i T_{x_i}$, and call
it a {\it triangular sum}.  For a vector $b$ of length $k$ we define the generating function
$$
F(q):=F_b(q):=  \ssum{x\in \Z^k}{} q^{t(x)}=\ssum{n=0}{\infty} \co_b(n) q^n,
$$
where $\co_b(n)$ is the number of solutions to $t(x)=n$.  
We will omit the subscript of $\co_b(n)$ when it is clear from the context.  We will furthermore use $r(n)$ to denote the number of representations of $n$ by the corresponding (diagonal) quadratic form $\sum b_i x_i^2$ and $r_o(n)$ to denote the number of those representations with all $x_i$ odd.
For ease of notation, we will denote the triangular sum corresponding to $b$ with $[b_1,b_2,\dots,b_k]$ and the corresponding quadratic form by $(b_1,\dots, b_k)$.  

The Hurwitz class number for the imaginary quadratic order of discriminant $D<0$ will play an important role in our analysis below.    We recall the definition here.  
For a negative discriminant $D$, the Hurwitz class number $H(D)$ is the weighted number of equivalence classes of, not necessarily primitive, positive definite binary quadratic forms of discriminant $D$, where the weights are $1$ except for classes of forms equivalent to a multiple of $\left(x^2+y^2\right)$, which are counted with weight $\frac{1}{2}$, and for classes of forms equivalent to a multiple of $\left(x^2+xy+y^2\right)$, which are counted with weight $\frac{1}{3}$.
Every quadratic form of discriminant $D$ is a multiple of a primitive form of
discriminant $D'=D/f^2$, and the weights are reciprocal to $w(D')/2$, half the number of units in the unique order of discriminant $D'$, or, accordingly, to half the number of representations of the integer $1$ by the primitive form.  The usual class numbers $h(D)$ are hence related to the Hurwitz class number by
$$
H(D) = \sum_{f^2\mid D} \frac{h\left(\frac{D}{f^2}\right)}{\frac{1}{2}w\left(\frac{D}{f^2}\right)}.
$$
For an integer $n$, we will set $a_n:=\frac{v_3(n+1)}{\log_3(n+1)}$, so that $3^{v_3(n+1)}= (n+1)^{a_n}$ gives the $3$-part of $n+1$ as a power of $n+1$.

\begin{theorem}\label{repsizethm}
For $\epsilon>0$, there is an absolute constant $c_{\epsilon}$ such that 
if the triangular sum $t(x)$ represents $1,2,4,5$, and $8$,
then $t(x)$ represents every nonnegative integer $n$ at least 
$\min \{ c_{\epsilon} n^{\frac{1}{2}-\epsilon}, 16 n^{1-a_n}\}$ times.  
In particular, if $n$ is sufficiently large and 
$a_n < \frac{1}{2}$ then $t(x)$ represents $n$ at least $c_{\epsilon} n^{\frac{1}{2}-\epsilon}$ times.
\end{theorem}
\begin{proof}
We proceed with \begin{it}escalator lattices\end{it} as in \cite{Bhargava1}. Without loss of generality we 
have $b_1\leq b_2\leq\dots\leq b_k$.  Fixing $b=[b_1,\dots,b_{k-1}]$, 
we will \begin{it}escalate\end{it} to $[b_1,\dots,b_k]$ by making all 
possible choices of $b_k\geq b_{k-1}$ for which it is possible to 
represent the next largest integer not already represented. 
We will then develop an \begin{it}escalator tree\end{it} by 
forming an edge between $b$ and $[b_1,\dots,b_k]$, with 
$\emptyset$ as the root.  If $\sum_ib_iT_{x_i}$ represents every integer, 
then $b$ will be a leaf of our tree.

Since $\co(1)>0$, it follows that $b_1=1$.  We need $\co(2)>0$, so $b_2=1$ or $b_2=2$.  If $b_2=1$, then we need $\co(5)>0$, so $1\leq b_3\leq 5$.  For $b_3=3$, we need $\co(8)>0$, so $3\leq b_4\leq 8$. Likewise, if $b_2=2$, then $2\leq b_3\leq 4$.  Therefore, if $\co(n)>0$ for every $n$, then we must have one of the above choices of $b_i$ as a sublattice.  By showing that each of these choices of $b_i$ satisfies $\co(n)>0$ for every $n$, we will see that this condition is both necessary and sufficient.  

All of the cases other than $[1,1,3,k]$ with $3\leq k\leq 8$ are covered by Liouville's Theorem.  However, to obtain the more precise version given in Theorem \ref{repsizethm}, we will use quadratic form genus theory.

One sees easily that 
$$
q^{\ssum{i=1}{k}b_i}F(q^8)=\ssum{x}{} q^{\ssum{i=1}{k} b_i (2x_i+1)^2 },
$$
so that $\co(n)=r_o\left(8n-\ssum{i=1}{k}b_i\right)$.  For the forms $b=[1,1,1]$, $[1,1,4]$, $[1,1,5]$, $[1,2,2]$, and $[1,2,4]$, congruence conditions modulo $8$ imply that 
$$
r_o\left(8n-\ssum{i=1}{k}b_i\right)=r\left(8n-\ssum{i=1}{k}b_i\right).
$$
Moreover, for each of these choices of $b$, $(b_1,b_2,b_3)$ is a genus 1 quadratic form.  Therefore, extending the classification of Jones \cite[Theorem 86]{Jones1} to primitive representations when the integer is not squarefree, $\co_{[1,1,1]}(n) = 24 H(-(8n+3))$, $s_{[1,1,4]}(n) = 4 H(-4(8n+6))$, $\co_{[1,2,2]}(n) = 4 H(-4(8n+5))$, and $\co_{[1,2,4]}(n) = 2 H(-8(8n+7))$.

For $[1,1,5]$ we must be slightly more careful since $5$ divides the discriminant.  We will explain in some detail how to deal with this complication and then will henceforth ignore this difficulty when it arises. 
 For $5\nmid 8n+7$ we have $\co_{[1,1,5]}(n) = 4 H(-5(8n+7))$.  Hence the only difficulty occurs with high divisibility by $5$.  For $p\neq 5$ the local densities are equal to those for bounded divisibility.  Thus, entirely analogously to the result of Jones we have $\co_{[1,1,5]}(n) = c_n H(-5(8n+7))$ for some constant $c_n>0$ which only depends $5$-adically on $8n+7$.  We calculate the cases $v_5(8n+7)\leq 3$ by hand.  Denote $5$-primitive representations of $m$ (i.e., $5\nmid \gcd(x,y,z)$) by $r^*(m)$.  Checking locally, for $5^2 \mid m:=8n+7$, we will obtain the result inductively by showing $\frac{r^*(25m)}{r^*(m)}=\frac{h(25m)/u(25m)}{h(m)/u(m)}$ and then summing to get $r(m)\geq  4 H(-5m)$.  
But, since $5\mid m$, we have $\frac{h(25m)/u(25m)}{h(m)/u(m)}=5$ by the class number formula (see \cite[Corollary 7.28, page 148]{Cox1}) so that this is a quick local check at the prime $5$.

Our proofs for $[1,1,2]$, $[1,2,3]$, and $[1,1,3]$ will be essentially identical.  For $[1,1,2]$, we note that if
$
x^2+y^2+2z^2=8n+4
$
has a solution with $x,y$, and $z$ not all odd, then taking each side modulo $8$ leads us to the conclusion that $x$, $y$, and $z$ must all be even.  Therefore, the solutions without $x$, $y$, and $z$ odd correspond to solutions of
$$
4x^2+4y^2+8z^2=8n+4,
\quad\textrm{that is, of}\quad
x^2+y^2+2z^2=2n+1.
$$
Using Siegel's theorem to compare the local density at $2$, we see that the average of the number of representations over the genus is three times as large for $8n+4$ as $2n+1$.  However, $(1,1,2)$ is again a genus $1$ quadratic form, so $r(8n+4)=3r(2n+1)$,
 and hence $\co_{[1,1,2]}(n)=r_o(8n+4)=r(8n+4)-r(2n+1)=2r(2n+1)$.  Thus by Theorem 86 of Jones \cite{Jones1} we have $\co_{[1,1,2]}(n) = 8H(-8(2n+1))$. 
 Similar arguments show that 
\begin{multline*}
\co_{[1,2,3]}(n)=r_{o,(1,2,3)}(8n+6)=r_{(1,2,3)}(8n+6)-r_{(4,2,12)}(8n+6)\\
=r_{(1,2,3)}(8n+6)-r_{(1,2,6)}(4n+3)=2r_{(1,2,6)}(4n+3).
\end{multline*}
Similarly  to the case $[1,1,5]$, we have $\co_{[1,2,3]}(n)\geq 2H(-12(4n+3))$.

For $[1,1,3]$ we see analogously that
$$
\co_{[1,1,3]}(n)=r_{o,(1,1,3)}(8n+5)=r_{(1,1,3)}(8n+5)-r_{(1,1,12)}(8n+5)=r_{(1,1,12)}(8n+5),
$$
and again $(1,1,12)$ is genus $1$.  
We conclude in the case $3\nmid (8n+5)$ that we have $\co_{[1,1,3]}(n)=4H(-3(8n+5))$, and we may henceforth assume that $3\mid 8n+5$ (i.e. $n\equiv 2\pmod{3}$). Local conditions imply that $3^{2j+1}(3\ell+2)$ is not represented by $(1,1,3)$, 
so we have escalated to $[1,1,3,k]$ for $k$ such that $3\leq k \leq 8$.  For $3\nmid k$, by choosing $x_4=1$ we have $\co_{[1,1,3,k]}(n)\geq 4H(-3(8(n-k)+5))$ since $3\nmid 8(n-k)+5$.
For $k=3$ we have 
$$
\co_{[1,1,3,3]}(n) = r_{(1,1,3,3)}(8(n+1)) +r_{(4,4,12,12)}(8(n+1)) - 2r_{(1,3,3,4)}(8(n+1)).
$$
Denoting the usual $d$-th degeneracy $V$-operator by $V(d)$ and the usual $U$-operator by $U(d)$ (cf. p. 28 of \cite{Ono1}), one may write the difference of the $\theta$-series $\sum_{n} r(8n) q^n$ for these quadratic forms as 
$$
\theta_{(1,1,3,3)}|U(8) + \theta_{(1,1,3,3)}|V(4)|U(8) - 2\theta_{(1,3,3,4)}|U(8).
$$ 
It is easy to conclude that the generating function $q F(z) = \sum_{n}\co_{[1,1,3,3]}(n)q^{n+1}$, with $q=e^{2\pi i z}$, is a weight 2 modular form of level $48$.  Using Sturm's bound \cite{Sturm1} and checking the first 16 coefficients reveals that $qF(z)=16\frac{\eta(2z)^4\eta(6z)^4}{\eta(z)^2\eta(3z)^2}$.  The coefficients are multiplicative, so that if we have the factorization $n+1=2^e 3^f \prod_{p>3} p^{e_p}$, then
$$
\co_{[1,1,3,3]}(n) = 2^{e+4} \prod_{p>3}\frac{p^{e_p+1}-1}{p-1}\geq 16\frac{n+1}{3^f} = 16(n+1)^{1-a_n}
$$

Finally, for $k=6$ we check $n<10$ by hand and then note that 
$$
\co_{[1,3,6]}(n) = r_{(1,3,6)}(8n+10)-r_{(2,3,6)}(4n+5),
$$
while both $(1,3,6)$ and $(2,3,6)$ are genus 1.  Hence for $n\not \equiv 2\pmod 3$ we have $\co_{[1,3,6]}(n)\geq 2H(-4(4n+5))$.  
We then take the remaining variable $x_4=1$ to obtain for $n\equiv 2\pmod{3}$ that $\co_{[1,1,3,6]}(n)\geq 2H(-4(4(n-1)+5))$, since $n-1\not\equiv 2\pmod{3}$.

Having seen that each of our choices of $b$ is indeed a leaf to the tree, we conclude that representing the integers $1,2,4,5,$ and $8$ suffices.
\end{proof}

\begin{remark}\rm
The constant $c_{\epsilon}$ in Theorem \ref{repsizethm} is ineffective because it relies on Siegel's lower bound for the class number, but the bound of $c_{\epsilon}n^{\frac{1}{2}-\epsilon}$ may be replaced with the minimum of finitely many choices of a constant times a Hurwitz class number of a certain imaginary quadratic order whose discriminant is linear in $n$.

We have the following example.  In this example, instead of considering $\co_{b}(n)$, we normalize the number of representations by
$$
\co_{b}'(n):=\frac{\co_b(n)}{2^k},
$$
where $k$ is the length of the sequence $b$.  This normalization is made so that $T_{x_i}=T_{-x_i-1}$ appears exactly once and in particular implies that $0$ is represented precisely once.  Using this normalization and the explicit bound in terms of the Hurwitz class number, we obtain for instance that if $1,2,4,5,$ and $8$ are represented, then the integer $195727301431$ is represented at least $270390$ times and the integer $48291403767737750$ is necessarily represented at least $90542761$ times (here $a_n\approx 0.364$), while the integer $50031545098999706=3^{35}-1$ is only necessarily represented once.  All of the bounds listed in these examples are sharp (i.e., there exists a triangular sum representing $1,2,4,5$, and $8$ which represents $195727301431$ precisely $270390$ times).
\end{remark}
\section{Cross Terms}\label{CrossTermsSection}
\noindent
Every quadratic polynomial $f$ in $k$ variables (over $\Q$) can be written uniquely as
$f(x)=Q(x)+\Lambda(x) + C$, where $Q(x)$ is a quadratic form in $k$ variables,
$\Lambda(x)$ is a linear form, and $C$ is a constant.  We will only consider 
quadratic polynomials such that $f(x)\in\integer$ for every $x\in \integer^k$.  The quadratic form $Q(x)$ is positive definite if and only if $f(x)$ is
bounded from below.  As in the introduction, $f(x_1,x_2,\dots, x_k)$ is a
normalized totally positive quadratic polynomial if $f$ is quadratic,
and the image of $\integer^k$ is contained in the non-negative integers 
while it contains $0$.
Clearly, for every positive definite quadratic form $Q(x)$ and linear form
$\Lambda(x)$ there is a unique $C\in\integer$ such that $f(x)=Q(x)+\Lambda(x) + C$ is normalized totally positive.

As noted before, $8T_x=(2x+1)^2-1=X^2-1$, if we put $X=2x+1$.
The
polynomial $X^2-1$ is normalized totally positive on the odd integers. 
With $Y=2y+1$, we find $8B_{xy}=4xy+2x+2y=XY-1$, where $B_{xy}:=\frac{1}{4}(2xy+x+y)$ is the
polynomial in $x, y$ satisfying $B_{xx}=T_x$. This way
$$8(aT_x+bT_y+cB_{xy})=aX^2+bY^2+cXY-(a+b+c).$$
If $C$ is the unique integer such that $aT_x+bT_y+cB_{xy}+C$ is
normalized totally positive, then $aX^2+bY^2+cXY+(8C-a-b-c)$ will be
the corresponding shifted quadratic form that is normalized totally
positive on the odd integers.

In order to describe our construction, we will say for simplicity that two quadratic polynomials $f_1$ and $f_2$ are \begin{it}(arithmetically) equivalent\end{it} if the number of solutions to $f_1(x)=n$ equals the number of solutions to $f_2(x)=n$ for every integer $n\geq 0$.

We will consider positive definite integral quadratic form (in $k$ variables) 
for which all cross terms in the matrix have {\it even} coefficients, 
so the cross terms of the quadratic form are $0\bmod 4$. 
This restriction is natural if
one keeps in mind that we are interested in the integers {\it oddly}
represented by forms. 

If $Q$ and $\widetilde{Q}$ are two equivalent quadratic forms such that the 
isomorphism preserves the condition that $X_i$ is odd, then we shall refer to 
them as \begin{it}equivalently odd\end{it}, and denote the equivalence class 
of such forms as $[Q]_{o}$.  

For any positive definite quadratic form with cross terms divisible by four,
we write
$$Q=a_1X_1^2+\cdots+a_kX_k^2+\sum_{i\neq j}4c_{ij}X_{i}X_{j},$$
we now define
$f_Q=f_{[Q]_{o}}$ to be the unique normalized totally positive quadratic 
polynomial
$$
f_Q:=a_1T_{x_1}+\cdots+a_kT_{x_k}+\sum_{i\neq j}4c_{ij}B_{x_ix_j}+C.
$$
We will refer to $f_Q$ as a \begin{it}triangular sum with cross terms\end{it}.

We will show that triangular sums with cross terms do not satisfy any 
finiteness theorem, and hence there is no overarching finiteness theorem for
quadratic polynomials, as stated in Proposition \ref{quadpolythm}.
To do so, for every positive integer $n$ we will 
construct a triangular sum with cross terms $f_n$ which represents 
precisely every non-negative integer other than $n$.

The following notation will be used. If $f$ and $g$ are polynomials
in $k$ and $\ell$ variables, we denote by $f\oplus g$ the sum of
the two as a polynomial in $k+\ell$ variables (so $f$ and $g$ are
assumed to share no variables).
\begin{theorem}\label{CrossCounterThm}
Let $\mathcal{Z}$ be a subset of the positive integers. For every proper subset 
$\mathcal{Y}\subsetneq \mathcal{Z}$ there exists a triangular sum with cross terms
representing $\mathcal{Y}$ but not representing $\mathcal{Z}$.
\end{theorem}

\begin{proof}
Let a proper subset $S_0$ of a given subset $S$
of the positive integers be given.  Choose a positive integer 
$n\in S\wout S_0$.  We will proceed by explicit construction of the 
triangular sum with cross terms $f_n$ which represents every integer 
other than $n$.  

First note that if the smallest positive integer {\it not} represented by 
$f$ is $n$, then, since the sum of three triangular numbers represents 
every non-negative integer, we have that 
$f\oplus (n+1)(T_x\oplus T_y\oplus T_z)$ 
represents all $m\not\equiv n\pmod{n+1}$.  But then we can choose 
$f_n:=f\oplus (n+1) (T_x\oplus T_y\oplus T_z) \oplus (2n+1) T_w$.  
It is therefore equivalent to construct $f$ for which $n$ is the smallest 
positive integer not represented by $f$.

Consider the quadratic form 
$$
Q^{(N)}(X,Y):=NX^2+NY^2+4XY,
$$
and denote the corresponding triangular sum with cross terms by  $f^{(N)}$; then
$$
f^{(N)}(x,y)=NT_{x}+NT_{y} + (2xy+x+y) +1.
$$
We first show that it is sufficient to determine that the generating function for $f^{(N)}$ is
\begin{equation}\label{Otermeqn}
2+2q+O(q^{N-12}).
\end{equation}
Assuming equation (\ref{Otermeqn}), then the generating function for 
$ g_n:=\oplus_{i=1}^{n} f^{(N)}$ is
$$
2^n\left(1+\binom{n}{1}q+\dots+\binom{n}{n}q^{n}\right) + O(q^{N-12}).
$$
If we choose $N>n+13$, then the
first integer not represented by $g$ is $n+1$.  Therefore, we can take
$f_n=g_{n-1}$; this also suffices for $n=1$ (if we interpret the empty direct sum $g_0$ as $0$).

We now show that the generating function satisfies (\ref{Otermeqn}).  
Note that
$f^{(N)}(0,-1)=f^{(N)}(-1,0)=0$, while $f^{(N)}(0,0)=f^{(N)}(-1,-1)=1$.
Now, without loss of generality, assume that $|x|\geq |y|$ and $x\notin \{0,-1\}$.  Then, 
$$
|2xy+x+y| \leq  2|x|^2+2|x| = 4 T_{|x|},
$$
so that 
$$
f^{(N)}(x,y)\geq NT_{x}-4T_{|x|} + NT_y.
$$
When $x\leq -2$ it is easy to check that $4T_{|x|}\leq 12 T_{x}$ so that 
$$
NT_{x}-4T_{|x|} \geq (N-12)T_{x}\geq N-12
$$
and when $x>0$ 
$$
NT_{x}-4T_{|x|} = (N-4)T_{x}\geq N-4,
$$
since $T_x\geq 1$ for $x\notin \{0,-1\}$.  Since $T_y\geq 0$, our assertion is verified.
\end{proof}
\noindent
It is important here to note how the above counterexamples differ from the proof when we only have diagonal terms, since this observation will lead us to the proof of Theorem \ref{CrossmBoundThm} when $m_{f}$ is bounded.  

Call a triangular sum with cross terms $f_Q$ (and also any corresponding $\widetilde{f_Q}$) a \begin{it}block\end{it} if the corresponding quadratic form $Q$ has an irreducible matrix.  We will build an escalator lattice by escalating (as a direct sum) by a block at each step.  In Section \ref{Thm8Section}, the breadth each time we escalated was finite, so that the overall tree was finite.  In the above proof, however, there were infinitely many inequivalent blocks which represent $1$, so that the breadth is infinite.  What was expressed in the above proof was that the supremum of these depths went to infinity as we chose $N$ increasing in terms of $n$ in the proof.

For 
$$
f(x) = \sum_{i=1}^k b_i T_{x_i} + \sum_{1\leq i < j\leq k} c_{ij}(2x_ix_j+x_i+x_j) + C
$$
we will say that $f$ has  \begin{it}(cross term) configuration\end{it} $c=(c_{ij})$.  Since the matrix of $f$ is irreducible and hence the corresponding adjacency matrix is connected, we can assume throughout (by a change of variables) that for each $j>1$ there exists $i<j$ with $c_{ij}\neq 0$.

\section{Bounded norm}
\noindent
We will now construct a natural norm on $f_{Q}$ such that restricting this norm will again give a finiteness result.  Let a positive definite quadratic form with even cross terms in the corresponding matrix,
\begin{equation}\label{QDef}
Q(x):=\ssum{i=1}{k}b_ix_i^2 + \ssum{i<j}{} 4c_{ij}x_ix_j
\end{equation}
be given.  We define
$$
\widetilde{f}(x):=\widetilde{f_Q}(x):=\ssum{i=1}{k} b_i T_{x_i} + \ssum{i<j, c_{ij}\geq 0}{} c_{ij}(2x_ix_j+x_i+x_j) + \ssum{i<j, c_{ij}< 0}{} c_{ij}(2x_ix_j+x_i+x_j+1).
$$
\begin{remark}\rm
Note that the constant $c_{ij}$ is added every time $c_{ij}<0$;
this may not seem canonical at first, but notice that if $Q'$ is the 
equivalent quadratic form obtained by replacing $x_1$ with $-x_1$, then we 
find that this choice leads to $\widetilde{f_Q}=\widetilde{f_{Q'}}$.
\end{remark}
We next define 
$$
\widetilde{m}_{\widetilde{f}}:=-\min_{x\in\integer^k} \widetilde{f}(x),
$$
which is added to obtain the unique (up to equivalence) normalized totally positive quadratic polynomial $f_Q=\widetilde{f_Q}+\widetilde{m}_{\widetilde{f}}$
corresponding to $Q$ .  Thus, we can define the norm 
$$
m_{f_Q}:=m_{[Q]_o}:=\min_{Q'\in [Q]_o} |\widetilde{m}_{\widetilde{f_{Q'}}}|.
$$
In a sense, this norm measures the distance between $f_Q$ and the closest $\widetilde{f_{Q'}}$ in the equivalence class, where the distance is merely given by the absolute value of the normalization factor required.  If $m_{f}$ is bounded, then we will again find that checking a finite subset will suffice.  We may now state the following more precise version of Theorem \ref{CrossmBoundThm}.
\begin{theorem}\label{thm}
Fix an integer $m$ and a subset $\mathcal{Z}$ of the positive integers.  Then there is a finite subset $\mathcal{Y}_{m}\subset \mathcal{Z}$ depending only on $m$ and $\mathcal{Z}$ such that every triangular sum with cross terms $f$ satisfying $m_{f}\leq m$ represents $\mathcal{Z}$ if and only if it represents $\mathcal{Y}_{m}$.

Moreover, for $\mathcal{Z}$ equal to the positive integers, we find $\max \mathcal{Y}_{m}\gg m^2$.
\end{theorem}
\begin{remark}\label{Crossmrmk}\rm
It may be of interest to investigate the growth of $\max \mathcal{Y}_m$ in terms of $m$ in the case where $\mathcal{Z}$ consists of all positive integers.  
The $m=0$ case is precisely Theorem \ref{thm8}.  Following the bounds given in the proof of Theorem \ref{CrossmBoundThm}, computational evidence suggests that $\mathcal{Y}_{1}(\integer_{>0})$
equals
\begin{eqnarray*}
\{1, 2, 3, 4, 5, 6, 8, 9, 10, 11, 12, 13, 14, 16, 17, 19, 20, 23, 24, 25, 26, 29, 32, 33, 34, 35, 38, 41,\\
46, 47, 48, 50, 53, 54, 58, 62, 63, 75, 86, 96, 101, 102, 113, 117, 129, 162, 195, 204, 233\}.
\end{eqnarray*}
\noindent
A proof of the above identity using the techniques of Bhargava and Hanke \cite{BhargavaHanke1} developed in the proof of the 290-Theorem may require a careful analysis of a possible Siegel zero.  To exhibit this difficulty, consider the sum $g(x,y,z)=T_x+2T_y+6T_z$.  In the construction of $\mathcal{Y}_{1}(\integer_{>0})$ the computations imply that there are infinitely many $Q$ with $m_{f_Q}=1$ for which $g\oplus f_Q$ represents every positive integer.  Hence we cannot merely check each case individually and must know information about the integers represented by $g$ independently.  

Although it seems that $g$ represents all odd integers, a proof of this appears to be beyond current techniques due to ineffective lower bounds for the class number (see \cite{Kane3}).  However, since a possible Siegel zero for $L(\chi_d,s)$ would give a lower bound for the class number when $d'\neq d$ (both fundamental), one may be able to show that $g$ represents at least one of $n$ or $n-1$ for every positive integer $n$, which would suffice for showing the above identity.
\end{remark}

We will first give an overview of the proof; details can be found in the next section. 

Fix a positive integer $m$.
As in the above remark, we will escalate with blocks.  We will first show that when $m_{f}\leq m$, the number of blocks that are not dimension $1$ in any branch of the escalator tree is bounded, and that there are only finitely many choices for the configuration of each block.  We will then proceed by defining $N(M_1,M_2,\dots,M_k,c)$
to be the smallest integer not represented by the totally positive quadratic polynomial corresponding to 
$$
\widetilde{f}(x):=\ssum{i=1}{k}M_i T_{x_i} + \ssum{i<j,c_{ij}\geq 0}{}c_{ij}(2x_ix_j+x_i+x_j)+ \ssum{i<j,c_{ij}< 0}{}c_{ij}(2x_ix_j+x_i+x_j+1).
$$
Our claim is then equivalent to showing that in the escalator tree 
$$
\sup_{M_1,\dots M_k,c} N(M_1,M_2,\dots,M_k,c)
$$
is finite.  To do so, we will effectively show that with the configurations of blocks of dimension greater than one fixed, the supremum with $M_i$ sufficiently large is finite and independent of the choice of $M_i$, and then fix $M_1\leq m_1$, and again show that the resulting supremum is independent of $M_2,\dots, M_k$, and so forth.  Since there are only finitely many such choices of $c$, the result comes from taking the maximum of each of these suprema.

\section{Proof}
\noindent
To prove Theorem \ref{thm}, and hence Theorem \ref{CrossmBoundThm}, we
begin with a lemma that will show that there are only finitely many choices of the cross term configuration.  
\begin{lemma}\label{clemma}
If $m_{f}\leq m$, then there are only finitely many choices of the cross term configurations $c_{ij}$ of all blocks of dimension greater than one, up to equivalent forms.
\end{lemma}

\begin{proof}
First note that $m_{f\oplus g}=m_{f}+m_{g}$, so that we can only have at most $m$ blocks $f$ with $m_{f}>0$, while we will see that $m_f>0$ unless $f$ is one dimensional (and hence the block is a constant times $T_x$).  It therefore sufficies to show that each block $f$ of dimension greater than one has $m_f>0$ and those with the restriction $m_f\leq m$ have bounded dimension and bounded coefficients in the configuration.  Fix the configuration $c$ of a block $\widetilde{f}$ with dimension $k$ such that $\widetilde{m}_{\widetilde{f}}=m_{f}$, namely a minimal element.  We will recursively show a particular choice of $x_i$ such that
$$
\widetilde{f}(x)\leq - \max \{ \max_{i,j} |c_{ij}|, k-1\},
$$
so that the max of the $c_{ij}$ is bounded by $m$, and the dimension is bounded by $m+1$.  

First set $x_1=0$.  Since $\widetilde{f}$ is a block, we know at step $j$ that there is some $i<j$ such that $c_{ij}\neq 0$.  Choose $i<j$ such that $|c_{ij}|$ is maximal.  If $x_i=0$, then we set $x_j=-1$ if $c_{ij}>0$ and $x_j=0$ otherwise.  If $x_i=-1$ then we set $x_j=0$ if $c_{ij}>0$ and $x_j=-1$ otherwise.  

Since all of our choices of $x_i$ are $0$ or $-1$ and $T_{-1}=T_{0}=0$, the integer represented is independent of the diagonal terms $M_i$.  Now we note that for $x_i,x_j\in \{0,-1\}$ we have $2x_ix_j+x_i+x_j=0$ if $x_i=x_j$ and $2x_ix_j+x_i+x_j=-1$ otherwise.  Therefore, if $x_i=x_j$, then from our definition of $\widetilde{f}$, the cross term corresponding to $c_{ij}$ adds $0$ if $c_{ij}\geq 0$ and adds $-|c_{ij}|$ otherwise.  If $x_i=0$ and $x_j=-1$, then the cross term adds $-|c_{ij}|$ if $c_{ij}\geq 0$ and adds $0$ otherwise.  Therefore by our construction above, we know that for $|c_{ij}|$ maximal, we have added $-|c_{ij}|$ to our sum, and we never add a positive integer, so the sum is at most $-|c_{ij}|$.  Moreover, since the block is connected, we have added at most $-1$ at each inductive step, so that the sum is at most $-(k-1)$.
\end{proof}
\noindent
For simplicity, in our escalator tree, we will ``push'' up all of the blocks to the top of the tree which are not dimension $1$.  To do so, we will first build the tree with all possible choices of blocks which are not dimension $1$, and then escalate with only dimension $1$ blocks from each of the nodes of the tree, including the root (the empty set).  Thus, every possible form will show up in our representation.  This tree (without the blocks of dimension $1$) is depth at most $m$ in the number of blocks, but is of infinite breadth.  Henceforth, we can consider the configuration $c$ to be fixed, and take the maximum over all choices of $c$.

We will now see that the subtree from each fixed node is of finite depth.  Consider the corresponding quadratic form $Q$.  First note that the generating function for $Q$ when all $x_i$ are odd is the generating function for $Q$ minus the generating function with some $x_i$ even, and the others arbitrary, which is simply the generating function for another quadratic form without any restrictions, taking $x_i\to 2x_i$.  Thus, we have the difference of $\theta$-series for finitely many quadratic forms, and hence the Fourier expansion is a modular form.  Now we simply note that any modular form can be decomposed into an Eisenstein series and a cusp form (cf. \cite{Ono1}).  Using the bounds of Tartakowsky \cite{Tartakowsky1} and Deligne \cite{Deligne1}, as long as the Eisenstein series is non-zero, the growth of the coefficients of the Eisenstein series can be shown to grow more quickly than the coefficients of the cusp form whenever the dimension is greater than or equal to $5$, other than finitely many congruences classes for which the coefficients of both the Eisenstein series and the cusp form are zero.  

Therefore, as long as the Eisenstein series is non-zero, there are only finitely many congruence classes and finitely many ``sporadic'' integers which are not represented by the quadratic form.  Thus, after dimension 5, there are only finitely many congruence classes and finitely many sporadic integers not represented by the form $f$.  If at any step of the escalation, any of the integers in these congruence classes is represented, then we have less congruence classes, and only finitely many more sporadic integers which are not represented, so that the resulting depth is bounded.  For the dimension 1 blocks, it is clear that the breadth of each escalation is finite, so there are only finitely many escalators coming from this node.  Therefore, it suffices to show that the Eisenstein series is non-zero.  

Again using Siegel's theorem \cite{Siegel2}, the coefficients of the Eisenstein series are simply a linear combination of the values given by the local densities of the quadratic forms from the above linear combination of $\theta$-series.  At every prime other than $p=2$, the local densities of the quadratic forms, of which we are taking the difference of $\theta$-series, are equal, so we only need to show that the difference of the local densities at $p=2$ is positive.  However, the difference of the number of local representations at a fixed $2$ power must be positive, since the integer is locally represented with $x_i$ odd, except possibly for finitely many congruence classes if a high $2$-power divides the discriminant.

Therefore, we can define $\widetilde{N}(M_1,\dots,M_k,c)$
to be the maximum of $N(M_1,\dots$, $M_k$, $M_{k+1}$, $\dots M_{l}$, $c)$, where $M_{k+1}$ to $M_{l}$ are the dimension $1$ blocks coming from the (finite) subtree of this node.  
We will show that $\widetilde{N}(M_1,\dots,M_k,c)$ is independent of the choice of $M_i$ whenever $M_i$ is sufficiently large by showing that the resulting subtrees are identical.  We need the following lemma to obtain this goal.  We will need some notation before we proceed.

For a set $T$, define 
$q^{T}:=\ssum{t\in T}{} q^t$,
a formal power series in $q$.
For fixed sets $S,T\subseteq\nature$, we will say that a form $f(x):=\sum b_i T_{x_i}$  \begin{it}represents $S/T$\end{it} if for every $s\in S$ the coefficient of $q^s$ in $q^{T}g(q)$ is positive, where $g(q)$ is the generating function for $f(x)$ given by $g(q):=\sum_{x\in \integer^k}q^{f(x)}$.

\begin{lemma}\label{SmodTLemma}
Let a (diagonal) triangular form $f$ be given.  Fix $S,T_1,T_2\subseteq \nature$ and $M\in \nature$ such that $\min_{n\in T_2}\geq M$.  Define $T:=T_1\cup T_2$.  Then there exists a bound $M_{T_1,S}$ and a finite subset $S_0\subseteq S$, depending only on $T_1$ and $S$ such that if $M>M_{T_1,S}$, then $f$ represents $S/T$ if and only if $f$ represents $S_0/T_1$.
\end{lemma}
\begin{proof}
We will escalate as in \cite{Bhargava1} with a slight deviation.  At each escalation node, there is a least element $s\in S$ such that $S/T_1$ is not represented by the form $f$ corresponding to this node.  As in \cite{Bhargava1}, we shall refer to $s$ as the \begin{it}truant\end{it} of $f$.  To represent $\{s\}/T_1$, we must have some $t_1\in T_1$ such that $s-t_1$ is represented by $f+bT_{x}$.  Therefore, for each $t_1<s$ we escalate with finitely many choices of $b$, and there are only finitely many choices of $t_1$.  Thus, the breadth at each escalation is finite, and our argument above using modular forms shows that the depth is also finite, so there are only finitely many choices of $s\in S$ which are truants in the escalation tree.  Take $S_0$ to be the set of truants in the escalation tree and define $M_{T_1,S}:=\max{s\in S_0} s+1$.  The argument above shows that representing $S/T_1$ is equivalent to representing $S_0/T_1$.  When following the above process with $T$ instead of $T_1$ whenever $M>M_{T_1,S}$, we will have the same subtree and the same truants at each step, so that representing $S/T$ is equivalent to representing $S/T_1$, and hence representing $S/T$ is equivalent to representing $S_0/T_1$.
\end{proof}
\begin{remark}\rm
It is of interest to note that if we replace ``(diagonal) triangular form'' with ``quadratic form'' (without the odd condition), then the proof follows verbatim, since the breadth is also finite, so that this can be considered a generalization of Bhargava's result that there is always a finite subset $S_0$ of $S$ such that the quadratic form represents $S$ if and only if it represents $S_0$, since this is obtained by taking $T_1=T=\{0\}$.
\end{remark}

Now consider $X_j:=\{x: x_i \text{ arbitrary for }i\leq j,\ x_i\in \{0,-1\}\text{ otherwise}\}$ and define
$T_{1,j}:=\{ f(x): x\in X_j\}$ and $T_{2,j}:=\{ f(x): x\notin X_j\}$.
We will use Lemma \ref{SmodTLemma} with $T_1=T_{1,j}$ and $T_2=T_{2,j}$ for each $0\leq j\leq k$.  To use the lemma effectively, we will show the following lemma.
\begin{lemma}\label{XjLemma}
There exist bounds $M_{X_j}^{(i)}$ depending only on $M_1,\dots, M_j,c$ such that if $M_i\geq M_{X_j}^{(i)}$ for every $i>j$, then the smallest element of $T_{2,j}$ is greater than $M_{T_1,\nature}$, where $M_{T_1,\nature}$ is as defined in lemma \ref{SmodTLemma}.
\end{lemma}
\begin{proof}
We will proceed by induction.  For $j=0$, we will take 
$
M_{X_0}^{(i)}= M_{T_{1,0},\nature} + 6\ssum{j}{}|c_{ij}|$.
Noting that for $|x_j|<|x_i|$ we have $\left|2\left(x_i-\frac{x_j}{2}\right)x_j\right|\leq x_i^2$, we get the inequality
$$
c_{ij}(2x_ix_j+x_i+x_j)\geq -|c_{ij}|(2 T_{|x_i|}+2T_{|x_j|}).
$$
The case $j=0$ then follows from the fact that for $x_i\notin\{0,-1\}$ we have $T_{|x_i|}\leq 3T_{x_i}$.

We now continue by induction on $j$.  For the corresponding quadratic form, we note that plugging in $x_1=\frac{-\ssum{j>1}{} c_{1j}x_j}{2M_1}$ gives the minimal value over the reals.  The quadratic form $Q'$ obtained by specializing this value of $x_1$ has rational coefficients with denominator dividing $2M_1$.  We therefore can consider $\widetilde{Q}:=4M_1\cdot Q'$, which is a quadratic form of the desired type.  Thus, we can use the inductive step for $\widetilde{Q}$.  But this gives a bound which minimizes $\widetilde{Q}$, and hence $Q'$, but an arbitrary choice of $x_1$ must give a value greater than or equal to this, so the result follows.
\end{proof}

Now, by our choice of $X_j$, $T_{1,j}$ is independent of $M_i$ for $i>j$, since $T_{x_i}=0$.  Thus, fix $c$ and take $M_i\geq M_{X_0}^{(i)}$.  Then the corresponding subtrees are independent of the choice of $M_i$, so that $\sup \widetilde{N}(M_1,\dots, M_k,c)$ is the unique largest truant in the subtree (effectively we may replace $M_i=\infty$).  We may now fix $M_1\leq M_{X_0}^{(1)}$, since there are only finitely many such choices.  With this $M_1$ fixed, we define $T_{1,1}$ as above, and again find bounds for the other $M_i$.  Continuing recursively gives the desired result, since we know that $k\leq m$, so there are only finitely many suprema that we take.

To show that $\max \mathcal{Y}_m(\integer_{>0})\gg m^2$, we consider again the construction of our counterexamples.  Consider $f(x,y):= \bigoplus_{i=1}^{m}f^{(N)}\oplus T_{y}$.  Since $T_r=\sum_{n=1}^{r} n$, for $N$ sufficiently large the smallest integer not represented by $f$ is clearly $T_{m+1}-1\gg m^2$.

\section*{Acknowledgements}
The authors would like to thank W.K. Chan for helpful comments as well as the anonymous referee for a detailed and helpful report.

\end{document}